\documentclass[12pt,leqno]{article} 
\usepackage{amsmath,amssymb,amscd}
\title{Presenting Schur algebras as quotients of the universal
enveloping algebra of $\gl_2$} 
\author{Stephen Doty and Anthony Giaquinto}
 
\newenvironment{pf}{{\em Proof.}}{\hfill$\square$\par\medskip}
 
\makeatletter 
\renewcommand{\subsection}{\@startsection{subsection}{2}{0mm}{\baselineskip}{-\fontdimen2\font}{\normalfont\normalsize\bfseries\setcounter{equation}{0}}}
\makeatother
 
\newenvironment{xlem}{{\bf Lemma\ }\em}{\par\medskip}
\newenvironment{xthm}{{\bf Theorem\ }\em}{\par\medskip}

\newcommand{\Z}{{\mathbb Z}}
\newcommand{\Q}{{\mathbb Q}}

\newcommand{\C}{{\mathbb C}}
\newcommand{\g}{{\mathfrak g}}

\newcommand{\End}{\operatorname{End}}

\newcommand{\Dist}{\operatorname{Dist}}

\newcommand{\GL}{{\sf GL}}
\newcommand{\SL}{{\sf SL}}
\newcommand{\gl}{\mathfrak{gl}}
\renewcommand{\sl}{\mathfrak{sl}}
\newcommand{\Ha}{{H_1}}
\newcommand{\Hb}{{H_2}}
\newcommand{\Hi}{{H_i}}
\newcommand{\divided}[2]{#1^{(#2)}}
\newcommand{\B}{\mathcal{B}}

\begin{document}
\maketitle
\begin{abstract} We give a presentation of the Schur algebras 
$S_\Q(2,d)$ by generators and relations, in fact a presentation which
is compatible with Serre's presentation of the universal enveloping
algebra of a simple Lie algebra.  In the process we find a new
basis for $S_\Q(2,d)$, a truncated form of the usual PBW basis.  We
also locate the integral Schur algebra within the presented algebra as
the analogue of Kostant's $\Z$-form, and show that it has an integral
basis which is a truncated version of Kostant's basis. 
\end{abstract}

\section{Introduction}
Consider a $2$-dimensional vector space $E$ over the rational field
$\Q$, which shall remain fixed throughout the paper, and fix a basis
$e_1, e_2$ for $E$. Via this basis for $E$ the group of linear
automorphisms of $E$ will be identified with the group $\GL_2(\Q)$ and
the vector space $\End(E)$ under Lie bracket $[x,y]=xy-yx$ will be
identified with the Lie algebra $\gl_2(\Q)$. The action of the group
on $E$ extends diagonally to an action on tensor space $E^{\otimes
d}$. The corresponding representation
$$
\sigma_d: \GL_2(\Q) \to \End(E^{\otimes d})
$$
extends by linearity to a representation 
$$
\sigma_d: \Q\GL_2(\Q) \to \End(E^{\otimes d})
$$
of the group algebra $\Q\GL_2(\Q)$, and the Schur algebra $S(2,d)$ is
precisely the image of the representation. By differentiating
$\sigma_d$ we obtain a representation $\rho_d$ of the Lie algebra
$\gl_2(\Q)$ which extends linearly to a representation 
$$
\rho_d: U(\gl_2) \to \End(E^{\otimes d})
$$
of the universal enveloping algebra, and its image is the same as the
image of the group algebra.  Thus $S(2,d)$ is a homomorphic image of
$U(\gl_2)$.

Now if the original representation $\sigma_d$ is restricted to the
subgroup $\SL_2(\Q)$ it is obvious that the image of its group algebra
in $\End(E^{\otimes d})$ will still be the same as the image of the
group algebra of $\GL_2(\Q)$ since the two groups differ only in
scalars and $\End(E^{\otimes d})$ already contains all the scalars. So
the restriction of $\rho_d$ gives a representation of $\sl_2$ and it
follows that $S(2,d)$ is also a homomorphic image of $U(\sl_2)$.

Serre has given a presentation of $U(\g)$ for any semisimple Lie
algebra $\g$ and, in particular, for $U(\sl_2)$.  This is easily
adapted to give a presentation of $U(\gl_2)$ (see \S\ref{B}). Thus the
natural question arises: to find an efficient set of generators of the
kernel of $\rho_d$ in terms of the Lie algebra generators of
$U(\sl_2)$ or $U(\gl_2)$. In other words, what additional relations
must be imposed on the generators of the universal enveloping algebra
in order to define a presentation of the Schur algebra?  This is the
same as the question of finding a set of generators for the
annihilator of the $U$-module $E^{\otimes d}$ ($U=U(\gl_2)$ or
$U(\sl_2)$).

In this paper we obtain a precise answer to this question. Write
$e_{ij}$ for the usual matrix units, defined in terms of Kronecker's
delta by $e_{ij}=(\delta_{il}\delta_{jl})_{i,j}$.  In the Lie algebra
$\gl_2$, set $e=e_{12}$, $f=e_{21}$, $\Hi=e_{ii}$ for $i=1,2$, and
$h=H_1-H_2$.  One can easily compute the eigenvalues of the images
under the above representation of $\Hi$ and $h$ in $\End(E^{\otimes
d})$ and thereby determine the minimal polynomial of those
endomorphisms.  Then our result is that the relation given by the
minimal polynomial is precisely the additional relation needed to
describe the desired presentation of $S(2,d)$.  In other words, the
kernel of the quotient map $U(\sl_2) \to S(2,d)$ is generated by the
minimal polynomial of the image of $h$.  Similarly, the kernel of the
quotient map $U(\gl_2) \to \End(E^{\otimes d})$ is generated by the
minimal polynomials of the images of $\Ha$ and $\Hb$, along with one
additional relation, $\Ha + \Hb =d$, which can be used to eliminate
one of the two generators $\Hi$ from the presentation.  In proving
these results we obtain along the way a new basis for $S_\Q(2,d)$,
which is a truncated form of the usual PBW basis of the universal
enveloping algebra.

The algebra $S(2,d)=S_\Q(2,d)$ contains a $\Z$-order $S_\Z(2,d)$, the
integral Schur algebra, and if $k$ is any field we have an isomophism
$S_k(2,d) \simeq S_\Z(2,d)\otimes_\Z k.$ (See \cite{Green}.)  We show
that the integral Schur algebra is precisely the analogue of Kostant's
$\Z$-form in $S_\Q(2,d)$, i.e.\ it is the subring generated by all
divided powers of $e$, $f$. Moreover, it has an integral basis which
is precisely a truncated version of the usual basis of the Kostant
$\Z$-form.  (This basis is closely related to the work of Richard
Green \cite{RGreen}.)

Although we obtain these results over the rational field $\Q$ all the
arguments can be carried out equally well over any field of
characteristic zero, so the same presentation holds over any such
field, and, in particular, over the complex field $\C$.

For many problems in Lie theory, the rank $1$ case turns out to be of
fundamental importance for the general case, and that is why we devote
an entire paper just to that case.  We expect to treat the general
case in a later paper, in which one cannot expect to find such
explicit reduction formulas as those given here.  The quantized
version of these results is addressed in \cite{DG:quantum}, where
different techniques of proof are developed.

\section{Statement of results}\label{A}
We now describe our results more precisely.

Our first main result describes the Schur algebra $S(2,d)$ over the
rational field $\Q$ in terms of generators and relations. The result
is pleasant: $S(2,d)$ has the same presentation as $U(\sl_2)$ with
just one additional relation. 

\subsection{}\label{Aa}
\begin{xthm} Over $\Q$, the Schur algebra $S(2,d)$ is isomorphic 
with the associative algebra (with 1) generated by $e$, $f$, $h$ 
subject to the relations:
\begin{align}
& he-eh=2e; \qquad ef-fe=h; \qquad hf-fh=-2f;\label{Aa:a}\\ 
& (h+d)(h+d-2)\cdots(h-d+2)(h-d)=0.\label{Aa:b}
\end{align}
Moreover, this algebra has a ``truncated PBW'' basis over $\Q$ consisting of all
$f^a h^b e^c$ such that $a+b+c \le d$.
\end{xthm}

\noindent Note that taking $d\to \infty$ in the above recovers the
usual presentation of $U(\sl_2)$.  Moreover, it follows from the above
that for each $d$ there is a surjective quotient mapping $S(2,d+2) \to
S(2,d)$ defined by mapping generators onto generators.

By a linear change of variable ($h=2\Ha-d$ or $h=d-2\Hb$) we obtain
from the above theorem two equivalent presentations, which are more
convenient for describing the integral Schur algebra.

\subsection{}\begin{xthm}\label{Ac}
 Over $\Q$, the Schur algebra $S(2,d)$ is isomorphic 
with the associative algebra (with 1) generated by $e$, $f$, $\Ha$ 
subject to the relations:
\begin{align}
& \Ha e - e\Ha = e; \qquad ef-fe=2\Ha-d; \qquad \Ha f - f\Ha=-f;\label{Ac:a}\\ 
& \Ha(\Ha-1)\cdots(\Ha-d)=0.\label{Ac:b}
\end{align}
Moreover, this algebra has a ``truncated PBW'' basis over $\Q$ consisting 
of all $e^a \Ha^b f^c$ such that $a+b+c \le d$.
\end{xthm}

\subsection{}\begin{xthm}\label{Ad}
 Over $\Q$, the Schur algebra $S(2,d)$ is isomorphic 
with the associative algebra (with 1) generated by $e$, $f$, $\Hb$
subject to the relations:
\begin{align}
& \Hb e - e\Hb = -e; \qquad ef-fe=d-2\Hb; \qquad \Hb f - f\Hb=f;\label{Ad:a}\\ 
& \Hb(\Hb-1)\cdots(\Hb-d)=0.\label{Ad:b}
\end{align}
Moreover, this algebra has a ``truncated PBW'' basis over $\Q$ consisting 
of all $f^a \Hb^b e^c$ such that $a+b+c \le d$.
\end{xthm}

Of course Theorems \ref{Ac} and \ref{Ad} are equivalent to one
another via the relation $\Ha + \Hb = d$.

Our next result constructs the integral Schur algebra $S_\Z(2,d)$ in
terms of the generators given above.

\subsection{}\begin{xthm}\label{Ae}
The integral Schur algebra $S_\Z(2,d)$ is isomorphic 
with the subalgebra of $S_\Q(2,d)$ generated by all divided powers
$$
 e^{(m)}:= e^m/m!; \qquad f^{(m)}:= f^m/m! .
$$
Moreover, this algebra has a ``truncated Kostant'' basis over $\Z$
consisting of all
\begin{equation}\label{Ae:a}
f^{(a)} \dbinom{\Hb}{b} e^{(c)} \qquad (a+b+c \le d)
\end{equation} 
and another such basis consisting of all 
\begin{equation}\label{Ae:b}
e^{(a)} \dbinom{\Ha}{b} f^{(c)} \qquad (a+b+c \le d).
\end{equation}
\end{xthm}

It would appear to be useful to express the product of two basis
elements as a $\Z$-linear combination of basis elements. The next
result, combined with standard commutation identities in the
enveloping algebra, enables one to do this.

\subsection{}\begin{xthm}\label{Af}
In $S_\Z(2,d)$ we have the following reduction 
formulas for any nonnegative integers $a$, $b$, $c$:
\begin{align}
\divided{f}{a} \binom{\Hb}{b} \divided{e}{c} &= \sum_{k=s}^{\min(a,c)}
(-1)^{k-s}\binom{k-1}{s-1} \binom{b+k}{k} \divided{f}{a-k} \binom{\Hb}{b+k}
\divided{e}{c-k}\label{Af:a}\\ 
\divided{e}{a} \binom{\Ha}{b} \divided{f}{c} &= \sum_{k=s}^{\min(a,c)}
(-1)^{k-s}\binom{k-1}{s-1} \binom{b+k}{k} \divided{e}{a-k} \binom{\Ha}{b+k}
\divided{f}{c-k}\label{Af:b}
\end{align}
where $s=a+b+c-d$.
\end{xthm}

\section{Enveloping algebras}\label{B}
In this section we fix some basic notation and quote some standard
results from the theory of enveloping algebras.

\subsection{}\label{Ba} 
We write 
\begin{equation*}
e=\begin{bmatrix}0&1\\0&0\end{bmatrix}\quad
f=\begin{bmatrix}0&0\\1&0\end{bmatrix}\quad
\Ha=\begin{bmatrix}1&0\\0&0\end{bmatrix}\quad
\Hb=\begin{bmatrix}0&0\\0&1\end{bmatrix}
\end{equation*}
for the canonical basis elements of the Lie algebra $\gl_2$. The
elements $e,f, h:=\Ha - \Hb$ form the canonical basis of the Lie
subalgebra $\sl_2$.

\subsection{}\label{Bb} 
The universal enveloping algebra $U(\gl_2)$ is the associative algebra
(with 1) generated by $e, f, \Ha, \Hb$ subject to the relations
\begin{align}
  \Ha \Hb &= \Hb \Ha, \label{Bb:a}\\
\Ha e - e \Ha = e, &\quad  \Ha f-f \Ha = -f,\label{Bb:b}\\
\Hb e - e \Hb = -e, &\quad \Hb f - f \Hb = f,\label{Bb:c}\\
  ef-fe &= \Ha - \Hb.\label{Bb:d}
\end{align}

Moreover, $U(\sl_2)$, the associative algebra (with 1) generated by
$e,f,h$ and satisfying the relations
\begin{equation} \label{Bb:e}
he-eh=e, \quad hf-fh=-f, \quad ef-fe=h, 
\end{equation}
is isomorphic with the subalgebra of $U(\gl_2)$ generated by $e, f,
h=\Ha-\Hb$. (This follows at once from the Poincare-Birkhoff-Witt
theorem.)

\subsection{}\label{Bc}  
From the Poincare-Birkhoff-Witt theorem it also 
follows that the algebra $U(\sl_2)$ (resp., $U(\gl_2)$) has a
$\Q$-basis (the PBW-basis) consisting of all
$$
f^a h^b e^c  \qquad\text{(resp., $f^a \Ha^{b_1} \Hb^{b_2} e^c$)}
$$ 
as $a,b,b_1,b_2,c$ range over the nonnegative integers.

\subsection{}\label{Bd}
For any nonnegative integer $m$ and any element $T$ of an associative
$\Q$-algebra with $1$ we define
$$
\divided{T}{m} = T^m/(m!), \quad \dbinom{T}{m} =
T(T-1)\cdots(T-m+1)/(m!)
$$
and we define these for negative $m$ to have the value $0$. Note that
the usual identity
\begin{equation}\label{Bd:a}
\binom{T}{m+1} = \binom{T-1}{m+1} + \binom{T-1}{m}
\end{equation}
holds true in this situation.

We consider the $\Z$-group $G_\Z=\SL_2$ (resp., $\GL_2$) and denote by
$\Dist(G_\Z)$ its algebra of distributions (see \cite[Part II,
(1.12)]{Jantzen}).  $\Dist(G_\Z)$ is the subring of
$U(\sl_2)$ (resp., $U(\gl_2)$) generated by all
\begin{equation}
\divided{f}{m}, \quad \divided{e}{m} \qquad 
\left(\text{resp., } \divided{f}{m},
\quad \divided{e}{m}, \quad \binom{\Ha}{m}, \quad \binom{\Hb}{m}\right).
\end{equation}
In the case $G_\Z = \SL_2$ the algebra of distributions coincides with
Kostant's $\Z$-form $U_\Z(\sl_2)$ in $U(\sl_2)$ and has a $\Z$-basis
consisting of all
\begin{equation}
\divided{f}{a} \binom{h}{b} \divided{e}{c}
\end{equation}
as $a,b,c$ range over the nonnegative integers.  In the case
$G_\Z=\GL_2$ the algebra of distributions has a $\Z$-basis consisting
of all
\begin{equation}
\divided{f}{a} \binom{\Ha}{b_1} \binom{\Hb}{b_2} \divided{e}{c}
\end{equation}
as $a,b_1,b_2,c$ range over the nonnegative integers.  We call
elements of this form {\em monomials} of degree $a+b_1+b_2+c$ and
height $a+c$.  Similar language will be applied to elements obtained
from elements of the above form by interchanging $\Ha$ with $\Hb$ and
$e$ with $f$.

The following lemma is well-known (see \cite{Kostant}) and will be
needed later on.

\subsection{}\begin{xlem} \label{Be} 
Given any polynomial $P(t)$ in $\Q[t]$ and any nonnegative integer
$k$, we have the following identities in $U(\gl_2)$:
\begin{gather}
e\, P(\Ha) = P(\Ha-1)\, e, \quad  e\, P(\Hb) = P(\Hb+1)\, e,\label{Be:a}\\
f\, P(\Ha) = P(\Ha+1)\, f, \quad  f\, P(\Hb) = P(\Hb-1)\, f,\label{Be:b}\\
e \divided{f}{k} =\divided{f}{k}e + \divided{f}{k-1}(\Ha-\Hb-k+1)\label{Be:c}\\
\divided{e}{k}f =f\divided{e}{k} + (\Ha-\Hb-k+1)\divided{e}{k-1}\label{Be:d}\\
f\divided{e}{k} = \divided{e}{k}f + \divided{e}{k-1}(\Hb-\Ha-k+1)\label{Be:e}\\
\divided{f}{k}e = e\divided{f}{k} + (\Hb-\Ha-k+1)\divided{f}{k-1}.\label{Be:f}
\end{gather}
\end{xlem}

\begin{pf} 
By linearity it suffices to prove the identities \eqref{Be:a} and
\eqref{Be:b} for $P$ of the form $t^m$ where $m$ is a nonnegative
integer.  But this follows for positive $m$ from
\ref{Bb}\eqref{Bb:b},\eqref{Bb:c} by induction on $m$, and for $m=0$
it is vacuous.  The next two identities follow from the relation
\ref{Bb}\eqref{Bb:d} by induction on $k$.  Finally, the last two
relations follow by symmetry, by interchanging $e$ with $f$ and $\Ha$
with $\Hb$.
\end{pf}

\section{The algebra $\B_d$}  
In this section we define an algebra $\B_d$ (over $\Q$) by generators
and relations. Eventually $\B_d$ will turn out to be isomorphic with
the Schur algebra $S(2,d)$.

\subsection{}\label{Ca} 
Fix a nonnegative integer $d$. We define $\B_d$ to be the algebra with
$1$ generated by $e,f,\Ha,\Hb$ subject to the relations
\ref{Bb}\eqref{Bb:a}--\eqref{Bb:d}, together with the additional
relations
\begin{gather}
\Ha + \Hb = d; \label{Ca:a} \\
0 = \Ha(\Ha-1)\cdots(\Ha-d). \label{Ca:b}
\end{gather}
Thus by definition $\B_d$ is a homomorphic image of $U(\gl_2)$. 

Note that in the presence of the relation \eqref{Ca:a} the
relation \eqref{Ca:b} can be replaced by the equivalent
relation
\begin{equation} \label{Ca:c}
0 = \Hb(\Hb-1)\cdots(\Hb-d).
\end{equation}

An important remark that will be useful in the sequel is that the
defining relations for $\B_d$ are invariant under the operation of
interchanging $e$ with $f$ and $\Ha$ with $\Hb$.  We shall refer to
this property as {\em symmetry}.  (This operation defines a Lie
algebra automorphism of $\gl_2$, which induces a corresponding
automorphism of $\B_d$.)

\subsection{}\label{Cb} 
It follows from the defining relations that in the 
algebra $\B_d$ we have
\begin{equation} \label{Cb:a}
\binom{\Ha}{d+1} = 0 = \binom{\Hb}{d+1}.
\end{equation}
More generally, we have
\begin{equation} \label{Cb:b}
0 = \binom{\Ha}{b_1} \binom{\Hb}{b_2} \qquad (b_1+b_2 \ge d+1). 
\end{equation}
One derives this first in the case $b_1+b_2=d+1$ from
\ref{Ca}\eqref{Ca:a} and \ref{Ca}\eqref{Ca:b}; the
statement for $b_1+b_2 > d+1$ then follows immediately.

\subsection{}\begin{xlem} \label{Cc}
In the algebra $\B_d$ we have 
$$
f^a\binom{\Hb}{b} = 0 = \binom{\Hb}{b}e^a, \quad 
e^a\binom{\Ha}{b} = 0 = \binom{\Ha}{b}f^a
$$
for any pair $a,b$ of nonnegative integers satisfying $a+b=d+1$. In
particular, $e^{d+1} = 0 = f^{d+1}$.
\end{xlem}

\begin{pf}
Set $b = d+1-a$ and proceed by induction on $a$. The case $a=0$ is
true by equations \ref{Cb}\eqref{Cb:a}. Supposing by
induction that $f^a\dbinom{\Hb}{b} = 0 = \dbinom{\Hb}{b}e^a$, we
obtain
\begin{align*}
0 &= f^a\binom{\Hb}{b} = f^a\binom{\Hb}{b}f \\ 
  &= f^{a+1}\binom{\Hb+1}{b} \quad(\text{by Lemma
  \ref{Be}})\\
  &= f^{a+1}\binom{\Hb}{b} + f^{a+1}\binom{\Hb}{b-1}\\
  &= f^{a+1}\binom{\Hb}{b-1}
\end{align*}
and
\begin{align*}
0 &= \binom{\Hb}{b}e^a = e\binom{\Hb}{b}e^a \\ 
  &= \binom{\Hb+1}{b} e^{a+1} \quad(\text{by Lemma
  \ref{Be}})\\
  &= \binom{\Hb}{b}e^{a+1} + \binom{\Hb}{b-1}e^{a+1}\\
  &= \binom{\Hb}{b-1}e^{a+1}
\end{align*}
and thus by induction the first two equalities are proved. The next
two equalities follow by symmetry, by interchanging $e$ with $f$ and
$\Ha$ with $\Hb$.  The last statement is an obvious special case of the
equalities which precede it.  The proof is complete.
\end{pf}

\subsection{}\begin{xlem} \label{Cd}
In the algebra $\B_d$ we have the identities
\begin{align}
(\Ha-\Hb)\binom{\Hb}{b} &= (d-2b)\binom{\Hb}{b} - (2b+2)\binom{\Hb}{b+1}\\
(\Hb-\Ha)\binom{\Ha}{b} &= (d-2b)\binom{\Ha}{b} - (2b+2)\binom{\Ha}{b+1}
\end{align}
for any nonnegative integer $b$.
\end{xlem}

\begin{pf}
By \ref{Ca}\eqref{Ca:a} we have 
\begin{align*}
(\Ha-\Hb)\binom{\Hb}{b} & = (d-2\Hb)\frac{\Hb(\Hb-1)\cdots(\Hb-b+1)}{b!}\\
  & = ((d-2b) + (2b-2\Hb))\frac{\Hb(\Hb-1)\cdots(\Hb-b+1)}{b!}\\
  & = (d-2b)\binom{\Hb}{b} - 2\frac{\Hb(\Hb-1)\cdots(\Hb-b+1)(\Hb-b)}{b!}\\
  & = (d-2b)\binom{\Hb}{b} - (2b+2)\binom{\Hb}{b+1}
\end{align*}
which proves the first identity. The second follows by symmetry. 
\end{pf}

Recall the definition of height and degree for monomials given at the
end of \ref{Bd}.

\subsection{}\begin{xthm}\label{Ce}
In the algebra $\B_d$ all monomials of the form $\divided{e}{a}
\dbinom{\Ha}{b} \divided{f}{c}$ (resp., $\divided{f}{a}
\dbinom{\Hb}{b} \divided{e}{c}$) of degree $d+1$ are expressible as
$\Q$-linear combinations of monomials of the same form but of strictly
smaller degree and height.
\end{xthm}

\begin{pf}
By symmetry it is enough to prove the claim for monomials of the form
$\divided{e}{a} \dbinom{\Ha}{b} \divided{f}{c}$.  For convenience,
write $M(a,b,c)$ for such a monomial.  

The proof proceeds by induction on height. The truth of the claim for
monomials of height zero is the content of \ref{Cb}\eqref{Cb:a}.

Consider an arbitrary monomial $M(a,b,c)$ of degree $d+1$ and height
$s\ge1$.  We first consider the case $a\ge c$. Then $a$ must be larger
than $0$, and by an application of \ref{Be}\eqref{Be:a} and
\ref{Bd}\eqref{Bd:a} we obtain the equalitites
\begin{align*}
\divided{e}{a-1} \dbinom{\Ha}{b+1}e\divided{f}{c} &= 
\divided{e}{a-1}e\dbinom{\Ha+1}{b+1}\divided{f}{c} \\
&= e\divided{e}{a-1}\dbinom{\Ha}{b+1}\divided{f}{c} 
+ a\divided{e}{a}\dbinom{\Ha}{b}\divided{f}{c} .
\end{align*}
By rearranging the above we find that
\begin{equation}\label{Ce:a}
a M(a,b,c) = 
\divided{e}{a-1}\dbinom{\Ha}{b+1}e\divided{f}{c} - e M(a-1,b+1,c).
\end{equation}
By induction we know that $M(a-1,b+1,c)$ is expressible as a
$\Q$-linear combination of monomials of degree $\le d$ and height $\le
s-2$, so $e M(a-1,b+1,c)$ is expressible as a $\Q$-linear combination
of monomials of degree $\le d+1$ and height $\le s-1$.  By induction
(applied to the degree $d+1$ terms) this reduces to a $\Q$-linear
combination of monomials of degree $\le d$ and height $\le s-1$.

Hence our claim for $M(a,b,c)$ will follow once we show that the term
$\divided{e}{a-1}\dbinom{\Ha}{b+1}e\divided{f}{c}$ can be expressed in
the desired form. But, by \ref{Be}\eqref{Be:f} we have
\begin{align*}
M&(a-1,b+1,c) e = \divided{e}{a-1} \dbinom{\Ha}{b+1} \divided{f}{c} e\\
&= \divided{e}{a-1} \dbinom{\Ha}{b+1}\left(e\divided{f}{c} 
   + (\Hb-\Ha-c+1)\divided{f}{c-1} \right)\\
&= \divided{e}{a-1} \dbinom{\Ha}{b+1} e\divided{f}{c} 
   + \divided{e}{a-1} \dbinom{\Ha}{b+1}(\Hb-\Ha-c+1)\divided{f}{c-1}.
\end{align*}
By rearranging this we obtain the equality
\begin{equation}\label{Ce:b}
\begin{aligned}
\divided{e}{a-1}&\dbinom{\Ha}{b+1}e\divided{f}{c} = \\
&M(a-1,b+1,c) e
- \divided{e}{a-1} \dbinom{\Ha}{b+1}(\Hb-\Ha-c+1)\divided{f}{c-1}.
\end{aligned}
\end{equation}
in which the second term on the right-hand-side is zero if $c=0$, by the
notational conventions introduced in \ref{Bd}.  By Lemma \ref{Cd} this
second term is expressible as a linear combination of monomials of
degree $\le d+1$ and height $\le s-2$, and by induction it can be
expressed as a $\Q$-linear combination of monomials of degree $\le d$
and height $\le s-2$.  

Finally, by induction we know that $M(a-1,b+1,c)$ is expressible as a
$\Q$-linear combination of monomials of degree $\le d+1$ and of height
$\le s-2$, and by induction applied to the terms of degree $d+1$ this
is reducible to a $\Q$-linear combination of monomials of degree $\le
d$ and of height $\le s-2$.  Thus it follows that $M(a-1,b+1,c) e$ is
expressible as a $\Q$-linear combination of monomials of degree $\le
d+1$ and of height $\le s-1$ (commuting the $e$ to the front does not
increase degree or height), so by applying induction one more time to
the terms of degree $d+1$ we are finished, in the case $a \ge c$.

The argument in the other case $c\ge a$ is entirely similar and is
omitted.
\end{pf}

\section{Identifications} 
In this section we will prove that $\B_d$ is isomorphic with
$S(2,d)$. Theorems \ref{Aa}, \ref{Ac}, \ref{Ad} will then follow.

\def\HHa{\overline{\Ha}}
\def\HHb{\overline{\Hb}}
\subsection{}\begin{xlem}\label{Da}
Write $\overline{X}$ for the image of $X\in U(\gl_2)$ under the
representation $\rho_d: U(\gl_2) \to \End(E^{\otimes d})$. Then 
we have the following identities
\begin{align}
&\HHa+\HHb = d\\
&\HHa(\HHa-1)\cdots(\HHa-d) = 0\\
&\HHb(\HHb-1)\cdots(\HHb-d) = 0
\end{align}
in the image $S(2,d)$ of $\rho_d$.  In fact, $T(T-1)\dots(T-d)$ is the
minimal polynomial of $\overline{\Hi}$ for $i=1,2$.
\end{xlem}

\begin{pf}
Set $W=E^{\otimes (d-1)}$ so that $E^{\otimes d} \simeq E \otimes W$.
We have the equality 
$$
\rho_d(\Hi)=1_E\otimes \rho_{d-1}(\Hi) + \rho(\Hi)\otimes 1_W.
$$
Now clearly $\rho(\Ha)+\rho(\Hb)=\rho(\Ha+\Hb)=\rho(1)=1$. Moreover,
the eigenvalues of $\rho(\Hi)$ are $0,1$ so its minimal polynomial is
$T(T-1)$. Thus the results hold when $d=1$.

For $d>1$ we have by the above and induction that
\begin{align*}
\rho_d(\Ha)+\rho_d(\Hb) &= 1_E\otimes(\rho_{d-1}(\Ha)+\rho_{d-1}(\Hb))
+ (\rho(\Ha) + \rho(\Hb))\otimes 1_W\\
&= 1_E\otimes(d-1)1_W + 1_E\otimes 1_W \\
&= (d-1) 1_{E\otimes W} + 1_{E\otimes W} = d 1_{E\otimes W}.
\end{align*}
This proves the first part. Moreover, by induction we may assume that
the eigenvalues of $\rho_{d-1}(\Hi)$ are $0,1,\dots,d-1$ (not counting
multiplicities).  Now suppose that $f$, $g$ are diagonalizable linear
operators on vector spaces $E$, $W$, resp.  Then the eigenvalues of
the linear operator $1\otimes g + f\otimes 1$ on $E\otimes W$ are all
of the form $\lambda+\mu$ where $\lambda$ is an eigenvalue of $E$ and
$\mu$ is an eigenvalue of $W$.  Applying this fact to $\rho_d(\Hi)$ we
see that its eigenvalues are $0,1,\dots,d$. The proof is complete.
\end{pf}

\noindent{\bf Remark.} A similar argument shows also that the minimal
polynomial of $\overline{h}$ is $(T+d)(T+d-2)\cdots(T-d+2)(T-d)$.

\subsection{}\begin{xthm}\label{Db}
Over $\Q$, the Schur algebra $S(2,d)$ is isomorphic with $\B_d$. 
The set of all $e^a \Ha^b f^c$ such that $a+b+c\le d$ is a $\Q$-basis,
as is the set of all $f^a \Hb^b e^c$ such that $a+b+c\le d$.
\end{xthm}

\begin{pf}
By the preceding lemma we see that the surjection $\rho_d:U(\gl_2)
\to S(2,d)$ factors through $\B_d$, giving the commutative diagram
\begin{equation}
\begin{CD}
U(\gl_2) @>>> S(2,d)\\
@VVV           @AAA\\
\B_d @= \B_d
\end{CD}
\end{equation}
in which all arrows are surjections. Since the set of all
$e^a\Ha^{b_1}\Hb^{b_2}f^c$ spans $U(\gl_2)$, it also spans $\B_d$. But
$\Hb = d-\Ha$, so it follows that the set of all $e^a\Ha^b f^c$
already spans.  Now by Theorem \ref{Ce} and it follows that the set of
all $e^a\Ha^b f^c$ satisfying the constraint $a+b+c \le d$ must span
the algebra $\B_d$.  By a similar argument we see that the same
statement holds for the set of all $f^a \Hb^b e^c$ such that $a+b+c\le
d$.  

In either case our spanning set is in one-to-one correspondence with
the set of all monomials in $4$ variables of total degree $d$ (set one
of the variables equal to $1$ to get the correspondence). It is well
known that this number is the dimension of $S(2,d)$ (see
\cite{Green}). Hence the map $\B_d \to S(2,d)$ must be an isomorphism
and the constrained spanning sets are bases. The theorem is proved.
\end{pf}

\noindent {\bf Remark.} We will henceforth simplify the notation of
Lemma \ref{Da}, writing simply $e,f,\Ha,\Hb$ for the images of these
elements in $S(2,d)$.

\subsection{}\label{Dc}
Theorem \ref{Ac} follows as an immediate corollary of the above
result, simply by using the relation $\Ha+\Hb=d$ to eliminate $\Hb$
from the presentation of $\B_d$.  Similarly, by eliminating $\Ha$ we
obtain Theorem \ref{Ad}. 

We now prove Theorem \ref{Aa}. Start from Theorem \ref{Ad} and set
$h=d-2\Hb$. So $\Hb=(d-h)/2$. By putting this into the relations
\ref{Ad}\eqref{Ad:a}, \ref{Ad}\eqref{Ad:b} and clearing denominators
and multiplying the second relation by $(-1)^d$ one arrives at the
desired relations \ref{Aa}\eqref{Aa:a} and \ref{Aa}\eqref{Aa:b}. Now
the set of all $f^a (h-d)^b e^c$ such that $a+b+c \le d$ is a
$\Q$-basis for the algebra. But by expanding $(h-d)^b$ by the binomial
theorem we see that every such basis element is expressible as a
$\Q$-linear combination of the elements $f^a h^b e^c$ satisfying
$a+b+c \le d$. Thus this latter set spans the algebra. But it has the
same cardinality as the dimension of the algebra, and so the proof is
complete.

\section{Integral reduction}

The purpose of this section is to prove Theorem \ref{Ae}. For this we
will need to know that the coefficients of the $\Q$-linear
combinations appearing in Theorem \ref{Ce} are actually integers. 
This is the content of Theorem \ref{Af}, which will be proved first. 

Note that our proof does not rely on the results in the preceding
section.  Thus we could use the results of this section instead of
Theorem \ref{Ce} to give another proof of Theorems \ref{Aa}, \ref{Ac},
and \ref{Ad}.

\subsection{}\begin{xthm} \label{Ea}
We have in $\B_d$ the equalities
\begin{align}
\divided{f}{a}\binom{\Hb}{b}\divided{e}{c} &= \sum_{k=1}^{\min(a,c)}
(-1)^{k-1} \binom{b+k}{k} \divided{f}{a-k}\binom{\Hb}{b+k}
\divided{e}{c-k}\\
\divided{e}{a}\binom{\Ha}{b}\divided{f}{c} &= \sum_{k=1}^{\min(a,c)}
(-1)^{k-1} \binom{b+k}{k} \divided{e}{a-k}\binom{\Ha}{b+k}
\divided{f}{c-k}
\end{align}
for all triples $a,b,c$ of nonnegative integers satisfying the
constraint $a+b+c=d+1$.
\end{xthm}

\begin{pf}
This is by a double induction on $a$ and $c$.  The case $a=c=0$ is
\ref{Cb}\eqref{Cb:a} and the cases $c=0$ ($a,b$
arbitrary) and $a=0$ ($b,c$ arbitrary) are already given in Lemma
\ref{Cc}.

By symmetry it suffices to prove just the first equality of the
theorem, which can be rewritten in the form
\begin{equation} \label{Ea:c}
0 = \sum_{k=0}^{\min(a,c)} (-1)^k \binom{b+k}{k} \divided{f}{a-k}
\binom{\Hb}{b+k} \divided{e}{c-k}
\end{equation}
By induction this holds for some fixed triple $a,b,c$ satisfying
$a+b+c=d+1$ and for all such triples whose first and last components
are no larger than $a$ and $c$.  We show how to derive the result in
the case $a+1,b-1,c$.  The case $a, b-1, c+1$ is entirely similar and
will be omitted.

The idea is to multiply \eqref{Ea:c} on the right by $f$ and then to
commute $f$ all the way to the left.  Using Lemma \ref{Be}
we obtain
\begin{equation*}
\begin{aligned}
0 & = \sum_{k=0}^{\min(a,c)} (-1)^k \binom{b+k}{k} \divided{f}{a-k}
\binom{\Hb}{b+k} \divided{e}{c-k} f \\ 
  & = \sum_{k=0}^{\min(a,c)} (-1)^k \binom{b+k}{k} \divided{f}{a-k}
\binom{\Hb}{b+k} \left(f\divided{e}{c-k} +
Q\divided{e}{c-k-1}\right)
\end{aligned}
\end{equation*}
where $Q = \Ha-\Hb-(c-k)+1$. Then by Lemma \ref{Be} 
once again we arrive at the equality
\begin{equation} \label{Ea:d}
\begin{aligned}
0 &= \sum_{k=0}^{\min(a,c)} (-1)^k \binom{b+k}{k} \divided{f}{a-k}f
\binom{\Hb+1}{b+k} \divided{e}{c-k}\\
& +  \sum_{k=0}^{\min(a,c)} (-1)^k 
\binom{b+k}{k} \divided{f}{a-k} \binom{\Hb}{b+k} (\Ha-\Hb-c+k+1) 
\divided{e}{c-k-1}.
\end{aligned}
\end{equation}

Now from Lemma \ref{Cd} it follows that 
\begin{multline}\label{Ea:e}
\binom{\Hb}{b+k}(\Ha-\Hb-c+k+1) = (a-k-b)\binom{\Hb}{b+k}\\ 
-(2b+2k+2)\binom{\Hb}{b+k+1}
\end{multline}
where we have made use of the equality
$(d-2(b+k))-(c-k-1)=(d+1-b-c)-k-b=a-k-b$.

Thus by putting \eqref{Ea:e} and the expansion 
$$
\binom{\Hb+1}{b+k} = \binom{\Hb}{b+k} + \binom{\Hb}{b+k-1}
$$
into \eqref{Ea:d} we obtain the equality
\begin{equation}\label{Ea:f}
\begin{aligned}
0 &= f \sum_{k=0}^{\min(a,c)} (-1)^k \binom{b+k}{k} \divided{f}{a-k}
\binom{\Hb}{b+k} \divided{e}{c-k}\\
& + \sum_{k=0}^{\min(a,c)} (-1)^k 
\binom{b+k}{k} f\divided{f}{a-k} \binom{\Hb}{b+k-1} \divided{e}{c-k}\\
& + \sum_{k=0}^{\min(a,c)} (-1)^k \binom{b+k}{k} (a-b-k) \divided{f}{a-k}
\binom{\Hb}{b+k} \divided{e}{c-k-1}\\
& + \sum_{k=0}^{\min(a,c)} (-1)^k 
\binom{b+k}{k}(-2b-2k-2)\divided{f}{a-k} \binom{\Hb}{b+k+1} \divided{e}{c-k-1}.
\end{aligned}
\end{equation}
Now by induction we claim that the first and fourth sums in the above
equality are zero. This is clear for the first sum.  To see this claim
for the fourth sum, note that
$$
\binom{b+k}{k}(-2b-2k-2) = -2(b+1)\binom{b+k+1}{k}
$$
and so the fourth sum in \eqref{Ea:f} is equal to 
$$
-2(b+1) \sum_{k=0}^{\min(a,c)} (-1)^k 
\binom{b+k+1}{k}\divided{f}{a-k}\binom{\Hb}{b+k+1} \divided{e}{c-1-k}
$$
which is a multiple of the $a,b+1,c-1$ case of the theorem, and hence
is zero by induction.  (If $\min(a,c)=c$ then the last term in the sum
is zero, so one can replace $\min(a,c)$ by $\min(a,c-1)$ without
changing the value of the sum.)

Thus the equality \eqref{Ea:f} reduces to 
\begin{equation*}
\begin{aligned}
0 &= \sum_{k=0}^{\min(a,c)} (-1)^k 
\binom{b+k}{k} (a+1-k)\divided{f}{a+1-k} \binom{\Hb}{b+k-1} \divided{e}{c-k}\\
& + \sum_{k=1}^{1+\min(a,c)} (-1)^{k-1} \binom{b+k-1}{k-1} (a-b-k+1) 
\divided{f}{a+1-k}\binom{\Hb}{b+k-1} \divided{e}{c-k}
\end{aligned}
\end{equation*}
where we have shifted the index of summation in the second
sum. Writing $m$ for $\min(a,c)$ we thus obtain
\begin{equation}\label{Ea:g}
\begin{aligned}
0 &= \sum_{k=1}^m (-1)^k R 
\divided{f}{a+1-k}\binom{\Hb}{b-1+k} \divided{e}{c-k}\\
  & + (a+1)\divided{f}{a+1}\binom{\Hb}{b-1}\divided{e}{c}\\
  & + (-1)^m \binom{b+m}{m}(a-b-m) \divided{f}{a-m}\binom{\Hb}{b+m}
      \divided{e}{c-1-m}
\end{aligned}
\end{equation}
where 
\begin{align*}
R & = \binom{b+k}{k}(a+1-k) - \binom{b+k-1}{k-1} (a-b-k+1)\\
  & = (a+1) \binom{b-1+k}{k}.
\end{align*}
Note that the last term in \eqref{Ea:g} is zero if
$m=c$. Otherwise $m=a<c$ and the term takes the form
\begin{multline*}
(-1)^{a}\binom{b+a}{a}(-b)\binom{\Hb}{b+a}\divided{e}{c-1-a}\\ 
= (-1)^{a+1}\, (a+1)\binom{b+a}{a+1}\binom{\Hb}{b+a}\divided{e}{c-1-a}
\end{multline*}
since 
\begin{align*}
(-1)^a \binom{b+a}{a}(-b) &= (-1)^{a+1} \binom{b+a}{b}b\\ &= (-1)^{a+1}
\binom{b+a}{b-1}(a+1)\\ &= (-1)^{a+1} \,(a+1) \binom{b+a}{a+1}.
\end{align*}
This shows that all the terms in equation \eqref{Ea:g} have a
common factor of $(a+1)$.  Putting these terms together and dividing
by $(a+1)$ we obtain the equality
$$
0 =  \sum_{k=0}^M (-1)^k 
\binom{b-1+k}{k}\divided{f}{a+1-k}\binom{\Hb}{b-1+k} \divided{e}{c-k}
$$
where $M=m=\min(a,c)=\min(a+1,c)$ in case $m=c$ and
$M=m+1=a+1=\min(a+1,c)$ otherwise.  In either case, $M$ is the minimum
of $a+1$ and $c$, and so we have obtained at last the desired equality
which completes the induction.
\end{pf}

\subsection{}\begin{xthm} \label{Eb}
Suppose that $a,b,c$ are nonnegative integers such that $s=a+b+c-d$ is
positive.  We have in $\B_d$ the equalities
\begin{align}
\divided{f}{a}\binom{\Hb}{b}\divided{e}{c} &= \sum_{k=s}^{\min(a,c)}
(-1)^{k-s} \binom{k-1}{s-1}\binom{b+k}{k} \divided{f}{a-k}\binom{\Hb}{b+k}
\divided{e}{c-k}\\
\divided{e}{a}\binom{\Ha}{b}\divided{f}{c} &= \sum_{k=s}^{\min(a,c)}
(-1)^{k-s} \binom{k-1}{s-1}\binom{b+k}{k} \divided{e}{a-k}\binom{\Ha}{b+k}
\divided{f}{c-k}.
\end{align}
\end{xthm}

\begin{pf}
By symmetry it is enough to prove the first equality. We proceed by
induction on $s$.  The case $s=1$ is precisely the content of Theorem
\ref{Ea}.  

Let $a,b,c$ be given such that $a+b+c-d=s+1$. If $a<s$ then $b+c-1\ge
d+1$, so by Lemma \ref{Cc} we have:
$$
\text{either }\quad \binom{\Hb}{b} \divided{e}{c-1} = 0 \quad\text{ or }\quad 
\binom{\Hb}{b-1} \divided{e}{c} = 0.
$$
In either case it follows that 
$$
\divided{f}{a}\binom{\Hb}{b}\divided{e}{c} = 0.
$$
By a similar argument one sees that this holds also if $c<s$.  Hence
we may assume that both $a$ and $c$ are $\ge s$. It is enough to argue
the result for the case $a\ge c\ge s$ since the other case $c \ge a
\ge s$ is similar.

Thus $a\ge 1$ and we have by the inductive hypothesis the equalities
\begin{align*}
\divided{f}{a}&\binom{\Hb}{b}\divided{e}{c} 
= \frac{f}{a} \divided{f}{a-1}\binom{\Hb}{b}\divided{e}{c} \\
&= \frac{f}{a} \sum_{k=s}^{\min(a-1,c)} (-1)^{k-s} \binom{k-1}{s-1} 
\binom{b+k}{k} \divided{f}{a-1-k} \binom{\Hb}{b+k} \divided{e}{c-k}\\
&= \frac{1}{a} \sum_{k=s}^{\min(a-1,c)} (-1)^{k-s} \binom{k-1}{s-1} 
\binom{b+k}{k}(a-k) \divided{f}{a-k} \binom{\Hb}{b+k} \divided{e}{c-k}\\
&= \frac{a-s}{a}\binom{b+s}{s} \divided{f}{a-s} \binom{\Hb}{b+s} 
\divided{e}{c-s} \\
&\qquad + \frac{1}{a} \sum_{k=s+1}^{\min(a-1,c)} (-1)^{k-s} \binom{k-1}{s-1} 
\binom{b+k}{k}(a-k) \divided{f}{a-k} \binom{\Hb}{b+k} \divided{e}{c-k}.
\end{align*}
Now the first term of the last equality above can be expanded by
Theorem \ref{Ea} (the base case of our present
induction) since $(a-s)+(b+s)+(c-s)=a+b+c-s=d+1$. Putting this in and
shifting the index of summation, we obtain the equalities
\begin{align*}
\divided{f}{a}&\binom{\Hb}{b}\divided{e}{c}\\ 
&= \frac{1}{a}\binom{b+s}{s}(a-s) \sum_{k=1}^{c-s} (-1)^{k-1}
  \binom{b+s+k}{k}\divided{f}{a-s-k}\binom{\Hb}{b+s+k}\divided{e}{c-s-k}\\
&\quad - \frac{1}{a} \sum_{k=s+1}^{\min(a-1,c)} (-1)^{k-(s+1)} 
\binom{k-1}{s-1}\binom{b+k}{k}(a-k) \divided{f}{a-k} \binom{\Hb}{b+k} 
\divided{e}{c-k}\\
&= \frac{1}{a}\binom{b+s}{s}(a-s) \sum_{k=s+1}^{c} (-1)^{k-(s+1)}
  \binom{b+k}{k-s}\divided{f}{a-k}\binom{\Hb}{b+k}\divided{e}{c-k}\\
&\quad - \frac{1}{a} \sum_{k=s+1}^{\min(a-1,c)} (-1)^{k-(s+1)} 
\binom{k-1}{s-1}\binom{b+k}{k}(a-k) \divided{f}{a-k} \binom{\Hb}{b+k} 
\divided{e}{c-k}.
\end{align*}
Now the second sum in the last equality above can be taken from $s+1$
to $c$ since $\min(a-1,c)$ is different from $c$ only if $a=c$, in
which case the additional term in the sum will be zero (the factor
$a-k$ is zero when $k=c=a$). Putting the two sums together and using
the computation
\begin{equation*}
\begin{split}
\frac{a-s}{a}&\binom{b+s}{s}\binom{b+k}{k-s} -
\frac{a-k}{a}\binom{k-1}{s-1}\binom{b+k}{k}\\
&= \frac{a-s}{a}\binom{k}{s}\binom{b+k}{k} -
\frac{a-k}{a}\binom{k-1}{s-1}\binom{b+k}{k}\\
&= \binom{k-1}{s}\binom{b+k}{k}
\end{split}
\end{equation*}
we obtain that
$$
\divided{f}{a}\binom{\Hb}{b}\divided{e}{c} = \sum_{k=s+1}^{c}
(-1)^{k-(s+1)} \binom{k-1}{s}\binom{b+k}{k} \divided{f}{a-k}\binom{\Hb}{b+k}
\divided{e}{c-k}
$$
and this completes the induction. 
\end{pf}

The following combinatorial result is well-known (see \cite{Kostant}).

\subsection{}\begin{xlem}\label{Ec}
Let $T_1,\dots, T_m$ be commuting indeterminates, and let
$F=F(T_1,\dots,T_m)$ be a polynomial over $\Q$ which takes on integer
values whenever the variables are replaced by integers. Then $F$ is
expressible as an integral linear combination of the polynomials 
$$
 \binom{T_1}{b_1} \binom{T_2}{b_2} \cdots \binom{T_m}{b_m}
$$
where $b_1, \dots, b_m$ all belong to the set of non-negative integers
and where $b_i$ does not exceed the degree of $F$ viewed as a
polynomial in $T_i$.
\end{xlem}

\subsection{}\label{Ed}
Theorem \ref{Eb} combines with the isomorphism $\B_d \simeq S(2,d)$
from Theorem \ref{Db} to prove Theorem \ref{Af}. 

We now prove Theorem \ref{Ae}. We begin with the representation 
\begin{equation}\label{Ed:a}
\rho_d: U(\gl_2) \to \End(E^{\otimes d}).
\end{equation}
We let $U_\Z(\gl_2)$ denote the subring of $U(\gl_2)$ generated by all
$\divided{f}{m}$, $\divided{e}{m}$, $\dbinom{\Ha}{m}$,
$\dbinom{\Hb}{m}$ for $m\ge0$. Then $U_\Z(\gl_2)$ is the same as the algebra
of distributions on the $\Z$-group $G_\Z=\GL_2$ from \ref{Bd}.

Now the map in \eqref{Ed:a} restricts to give a map $U_\Z(\gl_2) \to
\End(E^{\otimes d})$. Let $E_\Z$ be the $\Z$-submodule of $E$ 
spanned by the canonical basis elements $e_1$, $e_2$. It is easy to
see that $E_\Z$ is stable under the action of $U_\Z(\gl_2)$. It follows that
$E_\Z^{\otimes d}$ is also stable, and so the image of the above map
is contained within $\End(E_\Z^{\otimes d})$. Thus we have a
representation
\begin{equation}
\rho_d^\Z: U_\Z(\gl_2) \to \End(E_\Z^{\otimes d})
\end{equation}
By a result of Donkin \cite[2.1]{Donkin} (see also \cite{Doty:PRTCG})
the integral Schur algebra $S_\Z(2,d)$ is precisely the image of this
representation.

Hence $S_\Z(2,d)$ is spanned over $\Z$ by (the image of) the elements
$$
\divided{f}{a} \dbinom{\Ha}{b_1}\dbinom{\Hb}{b_2} \divided{e}{c}.
$$
But by Lemma \ref{Ec} we can express all $\dbinom{\Ha}{b_1}
\dbinom{\Hb}{b_2} = \dbinom{d-\Hb}{b_1}\dbinom{\Hb}{b_2}$ as integral
linear combinations of terms of the form $\dbinom{\Hb}{b}$. Thus
$S_\Z(2,d)$ is spanned by all $\divided{f}{a}\dbinom{\Hb}{b}
\divided{e}{c}$.  But by Theorem \ref{Eb} we see that the set of such
terms satisfying the restriction $a+b+c \le d$ is an integral spanning
set. But this set is linearly independent over $\Q$, hence linearly
independent over $\Z$, hence an integral basis. This proves
\ref{Ae}\eqref{Ae:a}.  Part (b) is proved similarly, starting with the
fact that the set of $\divided{e}{a}
\dbinom{\Ha}{b_1}\dbinom{\Hb}{b_2} \divided{f}{c}$ forms an integral
basis for $U_\Z(\gl_2)$.

Since $U_\Z(\sl_2)$ is generated by all divided powers
$\divided{e}{m}$, $\divided{f}{m}$ it is clear that $U_\Z(\sl_2)$ is
contained in $U_\Z(\gl_2)$.  It was shown in \cite{Donkin} (see also
\cite{Doty:PRTCG}) that the image of $U_\Z(\sl_2)$ in
$\End(E_\Z^{\otimes d})$ is the same as the image of
$U_\Z(\gl_2)$. This proves the first claim in Theorem \ref{Ae}. The
proof is complete.

\parskip=0pt\parindent=0pt\sf
\par Mathematical and Computer Sciences
\par Loyola University Chicago
\par Chicago, Illinois 60626 U.S.A.
\par E-mail: doty@math.luc.edu
\par \ \ \ \ \ \ \ \ \ tonyg@math.luc.edu

\end{document}